\documentclass[reqno]{amsart}
\usepackage{amsfonts, amssymb}
\usepackage{bm}
\usepackage{esint}

\newtheorem{theorem}{Theorem}
\newtheorem{lemma}{Lemma}

\newtheorem{proposition}{Proposition}
\theoremstyle{remark}
\newtheorem{remark}{Remark}
 
\newcommand{\field}[1]{\mathbb{#1}}  

\newcommand{\R}{\field{R}}

\def\DEF{\stackrel{\rm def}{=}}
\def\endproof{\hfill\hspace{-6pt}\rule[-14pt]{6pt}{6pt} \vskip22pt plus3pt minus 3pt}

\def\NN{\mathbb{N}}
\def\PP{\mathbb{P}}
\def\QQ{\mathbb{Q}}
\def\RR{\mathbb{R}}

\begin{document}

\parindent=0pt\parskip=8pt
\linespread{1.3}

\title[Discrete Entropy  of Generalized  Jacobi Polynomials] {Discrete Entropy  of Generalized  Jacobi Polynomials}

\author[A. Mart\'{\i}nez-Finkelshtein]{Andrei Mart\'{\i}nez-Finkelshtein}

\address[AMF]{Department of Mathematics, University of Almer\'{\i}a, Almer\'{\i}a, Spain}

\email{andrei@ual.es}

\author[P. Nevai]{Paul Nevai}

\address[PN]{Upper Arlington (Columbus), Ohio, USA}

\email{paul@nevai.us}

\author[A. Pe\~{n}a]{Ana Pe\~{n}a}

\address[AP]{Department of Mathematics and IUMA, University of Zaragoza,
Spain}

\email{anap@unizar.es}

\thanks{AMF and AP were supported by MICINN of Spain, jointly with the
European Regional Development Fund (ERDF), under grants
MTM2011-28952-C02-01 and MTM2012-36732-C03-02, respectively. Additionally, AMF was supported by
Junta de Andaluc\'{\i}a (Excellence Grant P11-FQM-7276 and the research group
FQM-229) and by Campus de Excelencia Internacional del Mar (CEIMAR) of the
University of Almer\'{\i}a. AP had partial support from DGA project E-64.
This work was completed during a visit of AMF to the Department of Mathematics of the Vanderbilt University. He acknowledges the hospitality of the hosting department, as well as a partial support of the Spanish Ministry of Education, Culture and Sports through the travel grant PRX14/00037.
}



\begin{abstract}

Given a sequence of orthonormal polynomials on $\RR$,$\{p_n\}_{n\geq 0}$, with $p_n$ of degree $n$, we define the discrete probability distribution $\bm \Psi_n(x) = \left(\Psi_{n,1}(x) , \dots \Psi_{n,n}(x) \right) $, 
 with $\Psi_{n,j}(x) = \big(\sum_{j=0}^{n-1} p_j^2(x)\big)^{-1}  p_{j-1}^2(x)$, 
 $j=1, \dots, n$.  In this paper, we study the asymptotic behavior as $n\to \infty$ of the Shannon entropy
$\displaystyle\mathcal{S} ((\bm \Psi_n(x))= -\sum_{j=1}^n  \Psi_{n,j}(x) \log ( \Psi_{n,j}(x))$, $x\in (-1,1)$, when the orthogonality weight is $ (1-x)^{\alpha}\, (1+x)^{\beta}\, h(x) $, $\alpha, \beta > -1$, and where $h$ is real, analytic, and positive on $[-1,1]$. We show that the limit 
$$
\lim_{n \to \infty} \left( \displaystyle\mathcal{S} ((\bm \Psi_n(x))- \log n\right)
$$
exists for all $x\in (-1,1)$, but its value depends on the rationality of $\arccos(x)/\pi$.

For the particular case of the Chebyshev polynomials of the first and second  kinds, we compare our asymptotic result with the explicit formulas for $\mathcal{S} (\bm \Psi_n(\zeta_j^{(n)}))$, where $\{ \zeta_j^{(n)}\}$ are the zeros of $p_n$, obtained previously in \cite{A-D-M-Y}.

\end{abstract} 

\date{\today}



\maketitle


\section {Introduction}

For a discrete probability distribution $\bm
\nu=(\nu_1,\nu_2,\cdots,\nu_n)$ with $\sum_{i=1}^n \nu_i=1$, we can
define its \emph{Shannon entropy} by
$$
\mathcal{S} (\bm \nu)= -\sum_{i=1}^n \nu_i \log (\nu_i), 
$$
that, by Jensen's inequality, satisfies
\begin{equation*}
0\le \mathcal{S}(\bm \nu)\le \log(n),
\end{equation*}
and the maximum of $\mathcal{S}(\bm \nu)$ is attained only at the uniform probability distribution
$$
\bm \nu^*=\left(\nu_1^*,\nu_2^*,\dots,\nu_n^* \right) = \left(1/n ,1/n ,\dots,1/n  \right) .
$$
Thus, along with the Shannon entropy, a natural measure of uncertainty
associated with a probability distribution $\bm \nu$ is its  ``\emph{distance}''
from $\bm \nu^*$, given by the \emph{directed} or \emph{Kullback--Leibler
divergence} 
\begin{equation}
\label{diver}
\mathcal{D} (\bm \nu)= \mathcal{D} (\bm \nu, \bm \nu^*)=\sum_{i=1}^n \nu_i \log \left( \frac{\nu_i}{\nu^*_i} \right)=  \log (n) - \mathcal{S}(\bm \nu) \geq 0. 
\end{equation}
 
Given a probability Borel measure $\mu$ supported on the real line $\RR$ with infinite number of points of increase and such that\footnote{We denote the set of nonnegative integers $\NN \cup \{0\}$ by $\NN_0$.}
$$
\int_\RR x^n d\mu(x) < \infty, \qquad n\in\NN_0,
$$
we can construct a sequence of orthonormal polynomials 
$$
p_n(x)=\kappa_n x^n+ \text{lower degree terms}, \qquad   \kappa_n >0, \quad n\in\NN_0,
$$
such that
\begin{equation*}
\int_\RR p_n(x)\,p_m(x)d\mu(x)=\delta_{n\,m}, \qquad n, m \in\NN_0 .
\end{equation*}
The corresponding reproducing kernel is
$$
K_n(x,y) = \sum_{j=0}^{n-1} p_j(x) p_j(y),
$$
that, for $x=y$, becomes the reciprocal of the $n$-th Christoffel function
\begin{equation*}
\lambda_n(x) \DEF \frac{1}{K_n(x,x)} \,.
\end{equation*}
For every $n\in \NN$ and $x\in \RR$ we can define the discrete probability distribution  
\begin{equation}
\label{probPsi}
	\bm \Psi_n(x) = \left(\Psi_{n,1}(x) , \dots \Psi_{n,n}(x) \right) , 
	\quad \text{with } \Psi_{n,j}(x) =  \lambda_n(x) p_{j-1}^2(x), 
	\quad j=1, \dots, n.
\end{equation}
Observe that this distribution does not depend on the normalization of the  measure $\mu$. 

\begin{remark}
Orthonormal polynomials $\{ p_n\}$ can be used also to define another sequence of probability distributions, $p_n^2(x) d\mu(x)$, defined on the support of the measure $\mu$.  The associated entropy,  
$$
- \int p_n^2(x) \log (p_n^2(x))\, d\mu(x),
$$
has been extensively studied, both its computation \cite{B-D-M04}, asymptotics \cite{A-D-M09} and applications in Physics, see e.g.~the survey~\cite{D-M-S01}.
\end{remark}

Our main goal is to study the asymptotic behavior of the Kullback--Leibler
divergence $\mathcal{D} (\bm \Psi_n(x))$, or, equivalently, that of the Shannon
entropy $\mathcal{S} (\bm \Psi_n(x))$, as $n\to \infty$ for $x$ in the bulk
of the support of the orthogonality measure $\mu$.  We restrict our attention
to absolutely continuous measures $\mu$ supported on a bounded interval of
$\R$, with $\mu'$ analytic and non-vanishing in the neighborhood of this
interval, except for the only possible singularities of a power type at the
endpoints of the support. Without loss of generality, we may assume that
\begin{equation}
	\label{mu}
	d\mu(x)=w(x)dx, \quad w(x) \DEF (1-x)^{\alpha}\, (1+x)^{\beta}\, h(x),
	\qquad x\in [-1,1],
\end{equation}
with $\alpha, \beta > -1$, and where $h$ is real, analytic, and positive on
$[-1,1]$. We call such kind of measures and the corresponding orthogonal
polynomials ``\emph{generalized Jacobi\/}''.

In what follows, when we have $x=\cos \theta \in (-1,1)$, then we also assume
that $\theta \in (0,\pi)$. 

One of the main results is the following theorem.

\begin{theorem}\label{thm1}
For $\mu$ given in \eqref{mu} and  $x=\cos \theta \in (-1,1)$, the limit
\begin{equation}
\label{mainlimit} 
	{\mathcal D}_\infty(x) \DEF \lim_{n\to \infty} \mathcal{D} (\bm \Psi_n(x))
\end{equation}
exists. Moreover,\footnote{Here and in what follows, $\operatorname{GCD}$
stands for the greatest common divisor.}
\begin{equation*}
	{\mathcal D}_\infty(x) = 
		\begin{cases}
		1- \log (2), & \text{ {\rm if} } \dfrac{\theta}{\pi} \notin \QQ , \\[3mm]
		\log (2) + 2\,  \widehat{\mathcal{S}}_{k,s},& \text{ {\rm if} }  \dfrac{\theta}{\pi}=\dfrac{s}{k} 
			\text{ {\rm with} } s, k \in \NN, \; s<k, \text{ {\rm and} } \operatorname{GCD} (s,k)=1 ,
		\end{cases}
\end{equation*}
where
\begin{equation}\label{S gorro}
	\widehat{\mathcal{S}}_{k,s} \DEF
	\displaystyle{\frac{1 }{k}\, \sum_{i=0}^{k-1} \, 
		\mathcal F\left( \cos\left((i+1/2)\frac{\pi s}{k}+ \varphi\left(\cos \frac{\pi s}{k} \right)-\pi/4\right)\right)},
\end{equation}
with
\begin{equation}\label{FF}
	\mathcal F(x) \DEF \begin{cases}
	x^2\log (x^2), & x>0, \\
	0, & x=0,
	\end{cases}
\end{equation}
and
\begin{equation} \label{psi}
	\varphi(x) \DEF
	\frac{1}{2}\big((\alpha+\beta)\theta-\alpha \pi\big)+
	\frac {\sqrt{1-x^2}}{2\pi} \fint_{-1}^1 \frac{\log (h(t))}{\sqrt{1-t^2}} \frac{dt}{t-x} \,.
\end{equation}

\end{theorem}

The integral in the right hand side of \eqref{psi} is understood in the sense of its principal value, that is,
$$
	\fint_{-1}^1 \frac{\log (h(t))}{\sqrt{1-t^2}} \frac{dt}{t-x} \DEF
	\lim_{\varepsilon \to 0} \left( 
	\int_{-1}^{x-\varepsilon} 
		\frac{\log (h(t))}{\sqrt{1-t^2}} \frac{dt}{t-x} + 
	\int_{x+\varepsilon}^1 \frac{\log (h(t))}{\sqrt{1-t^2}} \frac{dt}{t-x} 
	\right), \qquad x\in (-1,1).
$$

\smallskip

\begin{remark}
As formula~\eqref{entropia} below shows, the Shannon entropy $\mathcal{S} (\bm \Psi_n(x))$ (or  the Kullback--Leibler
divergence $\mathcal{D} (\bm \Psi_n(x))$) is closely related to the Christoffel function $\lambda_n(x)$ or the reproducing kernel $K_n(x,y)$. The latter exhibits a well-known universal behavior on the support of the orthogonality measure. In its most rudimentary form it is just the first limit in~\eqref{limite nucleo} below, while for the more sophisticated ``local'' version of this universality, leading to the sine kernel, see e.g.~ \cite{Lubinsky}. In all cases, the ``universal'' limit is continuous. This is no longer the case for the  Shannon entropy, as Theorem~\ref{thm1} illustrates, since function ${\mathcal D}_\infty$ is discontinuous everywhere in $(-1,1)$.
\end{remark}

\medskip

Let
$$
-1<\zeta_n^{(n)} <\dots < \zeta_1^{(n)}<1
$$
be the zeros of the $n$-th polynomial $p_n$.  In \cite{A-D-M-Y}, the authors studied the values of
$$
  \mathfrak{S}_{n,j}= \mathcal{S} (\bm \Psi_n(\zeta_j^{(n)})), \qquad j=1, \dots, n,
$$ 
finding explicit expressions for the case of orthonormal Chebyshev polynomials of the first and second kinds.
Recall that the orthonormal \emph{Chebyshev polynomials of the first kind} are given by the explicit formula
\begin{equation*}\label{explicitFirstKind}
   p_n(x)= T_n(x )=\begin{cases}  \dfrac{1}{\sqrt \pi}, & n=0\,, \\[3mm]
\sqrt{ \dfrac{2}{\pi}}\, \cos(n\theta), & n\in\NN,
\end{cases}\qquad x=\cos\theta\,,
\end{equation*}
for which $w(x)= ( 1-x^2 )^{-1/2}$ and
\begin{equation}\label{zeros_of_cheb}
\zeta _j^{(n)}=\cos \left(\frac{(2j-1)\pi }{2n}\right)\,,
\quad j=1, \dots, n\,,
\end{equation}
whereas the orthonormal \emph{Chebyshev polynomials of the second kind} are
\begin{equation*}\label{explicitSecondKind}
  p_n(x)=  U_n(x)= \sqrt{ \frac{2}{\pi}} \frac{\sin \left( (n+1) \arccos (x ) \right)}{\sqrt{1-x ^2}} =
 \sqrt{ \frac{2}{\pi}} \frac{\sin \left( (n+1) \theta \right)}{\sin (\theta )}\,, 
	\qquad x =\cos\theta\,, \quad n \in\NN_0\,,
\end{equation*}
with $w(x)= (1-x^2)^{1/2} $ and
\begin{equation}\label{zeros_of_2cheb}
\zeta_j^{(n)}=\cos \left( \frac{j \pi }{n+1}\right)\,,
\qquad j=1, \dots, n\,.
\end{equation}
Thus, it is interesting to study the compatibility of the results from
\cite{A-D-M-Y} with those stated in Theorem~\ref{thm1}. In other words, can
we reproduce \eqref{mainlimit}, ``stepping'' onto the zeros $\zeta_j^{(n)}$ 
only? The answer is yes, but not always.

Recall that the  explicit expression for the discrete entropy $\mathfrak{S}_{n,j}$ for orthonormal Chebyshev polynomials of the first kind  was derived in
\cite[Theorem~1, p.~99]{A-D-M-Y},
\begin{equation} \label{discreteentropy1}
	\mathfrak{S}_{n,j} =  \log n + \log 2 -1+\dfrac{\log 2}{n}-\mathcal R\left(\dfrac{d_n}{2n}\right) , 
\end{equation}
where $ d_n=\operatorname{GCD} (2j-1,n)$,
\begin{equation} \label{defRalt}
	\mathcal R(x) = - x \left( \psi \left(1-x\right) + 2 \gamma +
	\psi \left(1+x \right) \right) ,  
\end{equation}
$\gamma $ is the Euler-Mascheroni constant ($0.577\dots$), and
$\psi (x) \DEF \Gamma'(x)/\Gamma(x)$ is the digamma function. 
Alternatively, $\mathcal R$ can be evaluated using the series expansion, absolutely convergent for $ |x|<1$,
\begin{equation}\label{R-serie}
\mathcal R(x)=2 \sum_{k=1}^{\infty} \xi(2k+1) x^{2k+1},
\end{equation}
where $\xi(\cdot)$ is the Riemann zeta function.

An analogous expression was also obtained for the orthonormal Chebyshev polynomials of the second kind, see \cite[Theorem~2, p.~100]{A-D-M-Y}: 
\begin{equation}\label{discreteentropy2}
	\mathfrak{S}_{n,j} = \log (n+1) + \log 2 -1-\mathcal R\left(\frac{d_n}{n+1}\right),
\end{equation}
where now $ d_n=\operatorname{GCD} (j,n+1)$.

\begin{theorem}\label{teorema compatibilidad Ch-1}

Consider the orthonormal Chebyshev polynomials of the first or second kind and let 
$x=\cos \theta \in (-1,1)$. If $\theta/\pi \notin \QQ$, then there exists a
subsequence $\Lambda \subset \NN\times \NN$ such that
\begin{equation}
	\label{limitatzeros1}
	\lim_{(n,j)\in \Lambda} \zeta _j^{(n)} =x 
	\quad \text{\rm \&} \quad 
	\lim_{(n,j)\in \Lambda} \left(\mathfrak{S}_{n,j} - \mathcal{S}(\bm \Psi_n(x)) \right)=0.
\end{equation}
If  $\theta/ \pi= s/k$ where $s\in \NN$ and $k \in \NN$ with $s<k$  and 
$\operatorname{GCD} (s,k)=1$, then \eqref{limitatzeros1} still holds if the 
polynomials are of the second kind or if $k$ is even. However, for the
Chebyshev polynomials of the first kind and $k$ odd,  
\begin{equation}
	\label{weird}
	\limsup_{n} \left( \mathfrak{S}_{n,j_n} - \mathcal{S}(\bm \Psi_n(x)) \right)  <0
\end{equation}
for every subsequence  $\{j_n\} \subset \NN$.

\end{theorem}

\begin{remark}

A more precise statement than \eqref{weird} is given in \eqref{weird_strong}
below; it uses the function  $\mathcal R$ defined in \eqref{defRalt}--\eqref{R-serie}.

\end{remark}

\begin{remark}
The theorem above reveals a remarkable difference between the asymptotic behavior of the entropy of Chebyshev polynomials of the first and second kinds. A possible explanation is the fact that in the case of the polynomials of the first kind,  the denominator in the expression of the zeros \eqref{zeros_of_cheb} is always even, while for the second kind it can take any integer value, cf.~\eqref{zeros_of_2cheb}.
\end{remark}

Finally, we use the example of  Chebyshev polynomials of the first kind to compare the asymptotic values ${\mathcal D}_\infty(x) $ of the Kullback--Leibler
divergence for $x=\cos \theta \in (-1,1)$, when $\theta/\pi$ is either irrational (so that, according to Theorem~\ref{thm1}, ${\mathcal D}_\infty(x) = 1- \log (2)$), or rational. We see that ${\mathcal D}_\infty(x) $ attains neither its maximum nor its minimum at the irrational points:
\begin{proposition} \label{prop1}
Let $ \theta/\pi = s/k $, 
with $s, k \in \NN, \; s<k$, and $ \operatorname{GCD} (s,k)=1 $. Then for the orthonormal Chebyshev polynomials of the first kind,
$$
{\mathcal D}_\infty(x) = \begin{cases}
1- \log 2+\mathcal R\left(\dfrac{1}{k}\right) > 1- \log (2), & \text{if $k$ is even,} \\[3mm]
1- \log 2+ 2\left [\mathcal R\left(\dfrac{1}{2k}\right)-\dfrac{1}{2}\mathcal R\left(\dfrac{1}{k}\right) \right] < 1- \log (2), & \text{if $k$ is odd,}
\end{cases}
$$
where $\mathcal R$ is defined in \eqref{defRalt}--\eqref{R-serie}.
\end{proposition}

\section {Proof of Theorem \ref{thm1}}

Taking into account \eqref{diver} and \eqref{probPsi} we see that 
\begin{equation}\label{entropia}
\begin{split}
\mathcal{S} (\bm \Psi_n(x))&  = -  \log(\lambda_n(x))- \lambda_n(x) \sum_{i=0}^{n-1}\, p_{i}^2(x) \,\log(\, p_{i}^2(x)).
\end{split} 
\end{equation}

A crucial fact about the class of measures given in \eqref{mu} is that the
corresponding orthonormal polynomials satisfy the asymptotic
formula, valid uniformly on compact subsets of $(-1,1)$,
\begin{equation}\label{expresion polinomios}
	p_n(x)= 
	\sqrt{\frac{2}{\pi}} \frac{1}{\sqrt {w(x)}\, (1-x^2)^{1/4}} 
	\left(\cos \left((n+1/2)\theta+ \varphi(x)-\pi/4\right)+\mathcal O(1/n)\right), 
	\qquad x=\cos \theta,
\end{equation}
where the phase function $\varphi$ is given in \eqref{psi};
see \cite[(1.15) \& (1.33)]{K-M-V-V} where this asymptotics was proved using
the non-linear steepest descent method based on the Riemann--Hilbert
formulation of these polynomials.  

Given a generalized Jacobi $\mu$ as in \eqref{mu}, it is very well known that
\begin{equation}\label{limite nucleo}
	\lim_{n\to \infty} n\,  \lambda_n(x)=
	\pi\, w(x)\, \sqrt{1-x^2} 
	\quad \text{\&} \quad 
	\lim_{n\to \infty} \lambda_n(x)  p_n^2(x) =0 ,
\end{equation}
uniformly on compact sets of $(-1,1)$, see, e.g., \cite[Theorem 6.2.6 \&
Example~6.2.8, pp.~78--79]{N} for the first limit and \cite[Theorem~3.1.9,
p.~11]{N} or \cite[Theorem~2.1, p.~218]{NTZ} for the second one. Therefore,
we get from (\ref{entropia})--(\ref{limite nucleo}) that
\begin{equation}\label{entropia2}
	\mathcal{S} (\bm \Psi_n(x))
	= 
	\log \left(\frac{n}{2}\right) - \lambda_n(x) \sum_{i=0}^{n-1}\, p_{i}^2(x)   
	\log \left(\left( \cos \left( (i+1/2)\theta+ \varphi(x)-\pi/4 \right)+\epsilon_i(x) \right)^2\right)
	+
	o(1),
\end{equation}
where $\epsilon_i(x)=o(1)$ as $i\to \infty$ uniformly on compact sets of $(-1,1)$.

We have
\begin{equation}\label{descomposocion intermedia}
\begin{array}{ll}
	&\displaystyle{\lambda_n(x) \sum_{i=0}^{n-1}\, p_{i}^2(x)  
	\log \left(\left( \cos \big( (i+1/2)\theta+ \varphi(x)-\pi/4 \big)+\epsilon_i(x) \right)^2\right)}\\
	& = \displaystyle{\frac{2 \lambda_n(x)}{\pi w(x)\, \sqrt{1-x^2}} \sum_{i=0}^{n-1}\, \,
	\mathcal F\left(\cos \left( (i+1/2)\theta+ \varphi(x)-\pi/4 \right)+\epsilon_i(x)  \right)},
\end{array}
\end{equation}
where $\mathcal F$ is the function defined in \eqref{FF}.

Let us denote
$$
	y_i(x) \DEF \cos \left( (i+1/2)\theta+ \varphi(x)-\pi/4 \right), \qquad i=0, 1, \dots, n-1,
$$
and consider
\begin{align*}
	\sum_{i=0}^{n-1} \mathcal  F \left(y_i(x)+\epsilon_i(x)  \right) - \sum_{i=0}^{n-1} \mathcal  F \left(y_i(x)   \right) & = 
	\sum_{i=0}^{n-1} \mathcal  F' \left(y_i(x)+\nu_i(x) \epsilon_i(x)  \right) \epsilon_i(x) , 
	\qquad 0\leq \nu_i(x)\leq 1.
\end{align*}
Since $\mathcal F'$ is uniformly bounded on compact subsets of $[0,+\infty)$ and $\epsilon_i(x)=o(1)$, we can conclude that
\begin{equation}
\label{Lagrange}
\frac{1}{n}\sum_{i=0}^{n-1} \mathcal  F \left(y_i(x)+\epsilon_i(x)  \right) - \frac{1}{n}\sum_{i=0}^{n-1} \mathcal  F \left(y_i(x)   \right) = o(1), \quad n\to\infty,
\end{equation}
uniformly on compact subsets of $(-1,1)$.

Combining (\ref{entropia2})--\eqref{Lagrange}, we arrive at the asymptotic expression for the entropy
\begin{equation} \label{Sintermedio}
\mathcal{S} (\bm \Psi_n(x)) = \log \left(\frac{n}{2}\right)  -\frac{2 n \lambda_n(x)}{\pi w(x)\sqrt{1-x^2}}\, \widehat{\mathcal{S}}_n(x) +\, o(1),
\end{equation}
where
\begin{equation} \label{defShat}
\widehat{\mathcal{S}}_n(x) \DEF
 \frac{1 }{n}\, \sum_{i=0}^{n-1} \, \mathcal F\left(\cos \left( (i+1/2)\theta+ \varphi(x)-\pi/4 \right) \right).
\end{equation}

In order to prove Theorem~\ref{thm1}, we need to study the behavior of this function.

Assume first that $x=\cos \theta \in (-1,1)$ with $\theta/\pi\in \QQ$. Hence, there exist
$s \in \NN$ and $k \in \NN$ with $s<k$ and $\operatorname{GCD} (s,k)=1$, such that
$$
  \frac{\theta}{\pi}=\frac{s}{k} \, .
$$
Therefore, there exist non-negative integers $p$ and $q$ with $0\le q\le
k-1$ such that $n-1= p \,k +q$. We use the following straightforward lemma
(which is basically the idea behind the FFT algorithm), that can be proved by
direct calculation.

\begin{lemma}\label{lem1}

Let $g(n)$ be periodic with period $k$, that is, $g(n+k)=g(n)$ for all $n\in
\NN$. Let also $p$ and $q$, with $0\le q\le k-1$, be non-negative integers
such that $n-1= p \,k +q$. Then
$$
\frac{1}{n} \sum_{i=0}^{n-1} g(i)= \frac{1}{k} \sum_{i=0}^{k-1} g(i) + \frac{1}{n} \left(- \frac{q+1}{k}\sum_{i=0}^{k-1} g(i) + \sum_{i=0}^{q} g(i) \right).
$$
In particular, if $g(n)$ is uniformly bounded, then
\begin{equation*}
\lim_{n\to \infty} \frac{1}{n} \sum_{i=0}^{n-1} g(i)= \frac{1}{k} \sum_{i=0}^{k-1} g(i).
\end{equation*}

\end{lemma}

Applying Lemma~\ref{lem1} to
$$
g(n)=\mathcal F\left( \cos\big((n+1/2)\frac{\pi s}{k}+ \varphi(\cos \frac{\pi s}{k})-\pi/4\big)\right) ,
$$
we conclude that
\begin{equation*}
\lim_{n\to \infty} \widehat{\mathcal{S}}_n\left(\cos \frac{\pi s}{k}\right)= \frac{1 }{k}\, \sum_{i=0}^{k-1} \, \mathcal F\left( \cos \big((i+1/2)\frac{\pi s}{k}+ \varphi(\cos \frac{\pi s}{k})-\pi/4\big)\right) = \widehat{\mathcal{S}}_{k,s}. 
\end{equation*}
Together with \eqref{limite nucleo}  and \eqref{Sintermedio}, this
establishes the assertion of Theorem~\ref{thm1} for the case  $\theta/\pi \in
\QQ$.

Assume now that $x=\cos \theta \in (-1,1)$, but   $ \theta/\pi\notin \QQ$. 
By Kronecker's theorem (also known as Kronecker-Weyl's theorem), see, e.g.,
\cite[Theorem~IV, Chapter~III, p.~53]{Cassels}, the sequence\footnote{Here and in what follows, symbol $\lfloor\cdot\rfloor$ denotes  the mathematical integer part.}
$$
\left(n \frac{\theta}{\pi}- \left\lfloor n \frac{\theta}{\pi}\right\rfloor\right)_{n\in \NN},
$$
is dense and it is uniformly distributed in $[0,1]$. Thus, by \eqref{defShat},
\begin{equation*}
\lim_{n\to \infty} \widehat{\mathcal{S}}_n(\cos \theta)=  \int_0^1  \mathcal F\big( \cos ( \,y \pi +\frac{\theta}{2}+ \varphi(\cos \theta)-\pi/4)\big)dy.
\end{equation*}
Using the periodicity of the cosine function, we get
\begin{equation*}
\lim_{n\to \infty} \widehat{\mathcal{S}}_n(\cos \theta) = 
\int_0^1  \mathcal F\big( \cos  (\, y \pi )\big)dy=\frac{1}{2}-\log 2.
\end{equation*}
Again, combining this with \eqref{Sintermedio} we get the assertion of
Theorem~\ref{thm1} for the remaining case $x=\cos(\theta)$, $\theta/\pi
\notin \QQ$.
\endproof

\section{Proof of Theorem~\ref{teorema compatibilidad Ch-1} and Proposition~\ref{prop1}}

 \emph{Proof of Theorem~\ref{teorema compatibilidad Ch-1}.} 
 
 \textbf{1.} Let us start with the orthonormal Chebyshev polynomials of the first kind. Recall that $x=\cos \theta \in (-1,1)$.

\textbf{1.1.}
Assume first that $\theta/\pi\notin \QQ$. Using  Theorem~\ref{thm1}
and the explicit expression for  $\mathfrak{S}_{n,j}$ given in \eqref{discreteentropy1}--\eqref{R-serie},
we have
\begin{equation*} 
	\mathfrak{S}_{n,j} - \mathcal{S}(\bm \Psi_n(x)) = -\mathcal R\left(\frac{d_n}{2n}\right)+ o(1), 
	\qquad n\to \infty,
\end{equation*}
with $ d_n=\operatorname{GCD} (2j-1,n)$.

Denoting by $\PP\subset \NN$ the set of all prime numbers, we define  $\Lambda\equiv \Lambda_1
\subset \PP\times \NN \subset \NN\times \NN$ by
$$
	\Lambda_1 \DEF\left\{\left(p, \left\lfloor  \theta p/ \pi \right\rfloor \right) : \, p \in \PP \right\}.
$$
Clearly, we can make this sequence linearly ordered by agreeing that
$(p_1,j_1) \preceq (p_2, j_1) $ if $p_1\leq p_2$. Thus, when we take limits
along $(p,j)\in\Lambda_1$, we understand that $p\to \infty$.

By the construction, if $(p,j)\in \Lambda_1$, then 
$$
	\frac{\theta}{\pi}-\frac{3}{2p}\leq \frac{2 j  -1}{2 p} <\frac{\theta}{\pi}-\frac{1}{2p} \, ,
$$
so that
$$
	\lim_{ (p,j)\in \Lambda_1} \frac{2 j  -1}{2 p} = \frac{\theta}{\pi}<1 \, .
$$
Furthermore, $\operatorname{GCD}(2j-1,p)$  can take only two values, either
$1$ or $p$. Assuming $p>2$,  we must have $\operatorname{GCD}(2j-1,p)=p$ and, using that $2j-1<2p$, we
have $2j-1 = p$, so that $(p, (p+1)/2 )\in \Lambda_1$. If we assume that there
is an infinite subsequence of indices from $\Lambda_1$ of this form, we would get
$$
\lim_{ (p,j)\in \Lambda_1} \frac{2 j  -1}{2 p} = \lim_{  (p,j)\in \Lambda_1} \frac{p}{2 p} = \frac{1}{2} \, ,
$$
that contradicts our assumption that $ \theta/\pi\notin \QQ$. Thus,
for all sufficiently large $p$, we actually have $ \operatorname{GCD}(2j-1,p)=1$, and, therefore
\begin{equation*}
	\mathfrak{S}_{n,j} - \mathcal{S}(\bm \Psi_n(x)) = -\mathcal R\left(\frac{1}{2p}\right)+ o(1)= o(1), 
	\qquad n\to \infty,
\end{equation*}
that proves the assertion when $\theta/\pi\notin \QQ$.

\textbf{1.2.}
Assume now that $\theta/ \pi= s/k$ where $s\in \NN$ and $k\in \NN$ with $s<k$
and $\operatorname{GCD} (s,k)=1$. Using again Theorem~\ref{thm1} and  \eqref{discreteentropy1}--\eqref{R-serie} we have
\begin{equation} \label{difference3}
	\mathfrak{S}_{n,j} - \mathcal{S}(\bm \Psi_n(x)) = 
	2\log(2) -1 -\mathcal R\left(\frac{d_n}{2n}\right)+2 \,\widehat{\mathcal{S}}_{k,s}+ o(1), 
	\qquad n\to \infty,
\end{equation}
where we used the same notation as before.
From the explicit formula (\ref{S gorro}) for $\widehat{\mathcal{S}}_{k,s}$,
it is easy to see that in the case of the orthonormal Chebyshev polynomials of
the first kind,
\begin{equation} \label{equality-ks-k1}
\widehat{\mathcal{S}}_{k,s}=
 \frac{1 }{k}\, \sum_{i=1}^{k-1} \, \mathcal F\left( \cos \left(\frac{\pi i s}{k}\right)\right)= \frac{1 }{k}\, \sum_{i=1}^{k-1} \, \mathcal F\left( \cos \left(\frac{\pi i }{k}\right)\right) = \widehat{\mathcal{S}}_{k,1},
\end{equation}
where we have used that
\begin{equation}\label{reduccion s=1}
\left\{ \frac{ i s}{k} \mod 1:\, i= 1,2,\dots,k \right\}=\left\{ \frac{ i }{k} \mod 1:\, i= 1,2,\dots,k \right\}.
\end{equation}

\begin{remark}

In \cite[formula~(27), p.~108]{A-D-M-Y}, $\widehat{\mathcal{S}}_{n,1}$ is
normalized  in a different way 
because it lacks the normalizing factor $1/n$, so that our formulas will
slightly differ from those in \cite{A-D-M-Y}.

\end{remark}

\textbf{1.2.1.}
If $k$ is even, and, thus, $s$ is odd, we define $\Lambda\equiv\Lambda_2$ by
$$
	\Lambda_2 = \left\{ \left( \frac{k (2m+1)}{2}, \frac{s(2m+1)+1}{2} \right): m \in \NN \right\} \subset \NN\times \NN.
$$
Then
$$
	(n,j)\in \Lambda_2 \quad \Longrightarrow\quad \frac{2 j -1}{2 n}=\frac{s}{k} 
	\quad \text{\&}\quad 
	\frac{d_n}{2n}=\frac{1}{k} \, .
$$
Hence,
\begin{equation}  \label{asy1}
	\mathfrak{S}_{n,j} - \mathcal{S}(\bm \Psi_n(x)) = 
	2\log(2) -1 -\mathcal R\left(\frac{1}{k}\right)+2 \widehat{\mathcal{S}}_{k,1}+ o(1), 
	\qquad n\to \infty.
\end{equation}
Observe also that 
\begin{equation*}
	\widehat{\mathcal{S}}_{2m,1}=
 	\frac{1 }{2m}\, \sum_{i=1}^{2m-1} \, \mathcal F\left( \cos \left(\frac{\pi i }{2m}\right)\right)=  
	\frac{1 }{m}\, \sum_{i=1}^{m-1} \, \mathcal F\left( \cos \left(\frac{\pi i }{2m}\right)\right)
\end{equation*}
for $m\in \NN$, so that we can use formula (40) from \cite[Corollary~10, p.~111]{A-D-M-Y}, by which
\begin{equation} \label{identityR}
	2 \widehat{\mathcal{S}}_{k,1} = 1- 2\log(2) + \mathcal R\left(\frac{1}{k}\right)
\end{equation}
when $k$ is even, and \eqref{asy1} combined with \eqref{identityR} concludes the proof of
\eqref{limitatzeros1} for this case.

\textbf{1.2.2.}
Let us turn to the case when $k$ is odd. The key identity that holds in this
case is
\begin{equation} \label{newS}
	\widehat{\mathcal{S}}_{2m+1,1} = 
	\frac{1}{2}- \log(2) + \mathcal R\left(\frac{1}{2(2m+1)}\right) - \frac{1}{2}\mathcal R\left(\frac{1}{2m+1}\right),
	\qquad m\in \NN .
\end{equation}
Indeed,
\begin{align*}
\sum_{i=1}^{2m} \mathcal F\left( \cos \left(\frac{\pi i }{2(2m+1)}\right)\right) & =
	\sum_{i=1}^{m} \mathcal F\left( \cos \left(\frac{\pi i }{2m+1}\right)\right)+\sum_{i=1}^{m} \mathcal F\left( \cos \left(\frac{\pi (2i-1) }{2(2m+1)}\right)\right) \\
& =  \frac{1}{2} \sum_{i=1}^{2m} \mathcal F\left( \cos \left(\frac{\pi i }{2m+1}\right)\right)+\sum_{i=1}^{m} \mathcal F\left( \sin \left(\frac{\pi i }{2m+1}\right)\right) \\
& = \frac{2m+1}{2} \,  \widehat{\mathcal{S}}_{2m+1,1}+\frac{1}{2}\, \sum_{i=1}^{2m} \mathcal F\left( \sin \left(\frac{\pi i }{2m+1}\right)\right),
\end{align*}
so that
\begin{equation}
\label{identityS1}
\widehat{\mathcal{S}}_{2m+1,1} = \frac{2}{2m+1} \,\sum_{i=1}^{2m} \mathcal F\left( \cos \left(\frac{\pi i }{2(2m+1)}\right)\right) - \frac{1}{2m+1} \,\sum_{i=1}^{2m} \mathcal F\left( \sin \left(\frac{\pi i }{2m+1}\right)\right).
\end{equation}
By \eqref{identityR}, the first term in \eqref{identityS1} is $1-2\log(2) +
\mathcal R\left(\frac{1}{2(2m+1)}\right)$, while, by \cite[Proposition~13,
p.~114]{A-D-M-Y}, we have the following identity for the second term    
\begin{equation} \label{senoimpar}
 \frac{1}{2m+1} \,\sum_{i=1}^{2m} \mathcal F\left( \sin \left(\frac{\pi i }{2m+1}\right)\right)=\frac{1}{2}\left( 1-2\log(2) + \mathcal R\left(\frac{1}{2m+1}\right)\right), 
\end{equation}
that yields \eqref{newS}.

By \eqref{difference3} and \eqref{newS}, for $k$ odd, 
\begin{equation*} 
\begin{split}
\mathfrak{S}_{n,j} - \mathcal{S}(\bm \Psi_n(x)) &= 2\log(2) -1 -\mathcal R\left(\frac{d_n}{2n}\right)+2 \,\widehat{\mathcal{S}}_{k,1}+ o(1) \\
&= 2 \mathcal R\left(\frac{1}{2k}\right) -  \mathcal R\left(\frac{1}{k}\right) -\mathcal R\left(\frac{d_n}{2n}\right) + o(1), 
	\qquad n\to \infty,
\end{split}
\end{equation*}
where $ d_n=\operatorname{GCD} (2j-1,n)$. 

Observe that the coefficients in the power series expansions \eqref{R-serie} are
all positive, so that  $\mathcal R$ is convex on  $(0,1)$. As a
consequence, $\mathcal R(x/2)-\mathcal R(x)/2$ is decreasing  on $(0,1)$
and, therefore,
\begin{equation}
\label{convex function}
\mathcal R\left(\frac{x}{2}\right)-\frac{1}{2}\mathcal R(x)<0 \quad \text{for }x \in (0,1).
\end{equation}
In particular, for every choice of $\Lambda$,
\begin{equation}
	\label{weird_strong}
	\limsup_{(n,j)\in \Lambda}\left( \mathfrak{S}_{n,j} - \mathcal{S}(\bm \Psi_n(x)) \right) \leq 
	2 \mathcal R\left(\frac{1}{2k}\right) -  \mathcal R\left(\frac{1}{k}\right)<0,
\end{equation}
that establishes \eqref{weird}.

\textbf{2.}
Now we switch to the orthonormal Chebyshev polynomials of the second kind.
Let again $x=\cos \theta \in (-1,1)$.

\textbf{2.1.}
Assume $\theta/\pi \notin \QQ$. Using 
Theorem~\ref{thm1} and the explicit expression for  $\mathfrak{S}_{n,j}$
given in \eqref{defRalt}--\eqref{discreteentropy2}, we have
\begin{equation*}
	\mathfrak{S}_{n,j} - \mathcal{S}(\bm \Psi_n(x)) = 
	\log\left( \frac{n+1}{n}\right) -\mathcal R\left(\frac{d_n}{n+1}\right)+ o(1), 
	\qquad n\to \infty,
\end{equation*}
where $ d_n=\operatorname{GCD} (j,n+1)$. 

We build $\Lambda\equiv \Lambda_3 \subset \NN\times \NN$ and define linear ordering on it
similarly as it was done before for $\Lambda_1$. Namely, for each prime number $p$ take  $ j=
\left\lfloor  \theta (p-1)/ \pi \right\rfloor$  and then we denote all the
resulting pairs $(p-1,j)$ by $\Lambda_3$. By the construction,
$$
	\lim_{(p,j)\in \Lambda_3} \frac{  j  }{ p} = \frac{\theta}{\pi}<1,
$$
and $\operatorname{GCD} (j,p )=1$. Thus,
\begin{equation*}
	\mathfrak{S}_{n,j} - \mathcal{S}(\bm \Psi_n(x)) = 
	\log\left( \frac{p}{p-1}\right) -\mathcal R\left(\frac{1}{p}\right)+ o(1)= o(1), 
	\qquad n\to \infty.
\end{equation*}

\medskip

\textbf{2.2.}
Assume now that $\theta/ \pi= s/k$ where $s\in \NN$ and $k\in \NN$ with
$\operatorname{GCD} (s,k)=1$. Using  again Theorem~\ref{thm1} and the
explicit expression for  $\mathfrak{S}_{n,j}$ from \eqref{defRalt}--\eqref{discreteentropy2}, we obtain
\begin{equation*}
	\mathfrak{S}_{n,j} - \mathcal{S}(\bm \Psi_n(x)) = 
	\log\left( \frac{n+1}{n}\right) -1+2 \log(2)  -\mathcal R\left(\frac{d_n}{n+1}\right)+2\, \widehat{\mathcal S}_{k,s}+ o(1), 
	\qquad n\to \infty,
\end{equation*}
where $d_n=\operatorname{GCD} (j,n+1)$. In the case of  the orthonormal
Chebyshev polynomials of the second kind, $\widehat{\mathcal S}_{k,s}$
defined in  (\ref{S gorro}), has  the form
\begin{equation*}
	\widehat{\mathcal{S}}_{k,s}=
	\frac{1 }{k}\, \sum_{i=1}^{k} \, \mathcal F\left( \cos \left(\frac{\pi  s i}{k}-\frac{\pi}{2}\right)\right) = 
	\frac{1 }{k}\, \sum_{i=1}^{k-1} \, \mathcal F\left( \sin \left( \frac{\pi i}{k}\right)\right)=
	\widehat{\mathcal{S}}_{k,1},
\end{equation*}
where we have used \eqref{reduccion s=1} for the second equality.

If we take
$$
	\Lambda \equiv \Lambda_4 \DEF \left\{ \left( mk-1, sm \right): m \in \NN \right\} 
	\subset \NN\times \NN,
$$
then
\begin{equation*}
	(n,j)\in \Lambda_4 \quad \Longrightarrow\quad  \frac{ j}{ n+1}=\frac{s}{k} 
	\quad \text{\&}\quad 
	\frac{d_n}{n+1}=\frac{1}{k}
\end{equation*}
so that
\begin{equation} \label{difference11}
	\mathfrak{S}_{n,j} - \mathcal{S}(\bm \Psi_n(x)) = 
	2\log(2) -1 -\mathcal R\left(\frac{1}{k}\right)+2 \widehat{\mathcal{S}}_{k,1}+ o(1), 
	\qquad n\to \infty.
\end{equation}

\textbf{2.2.1} 
Let $k$ be even, say, $k=2m$. Observe  that
\begin{equation*}
	\widehat{\mathcal{S}}_{2m,1}=
	 \frac{1 }{2m}\, \sum_{i=1}^{2m-1} \, \mathcal F\left( \sin \left(\frac{\pi i  }{2m}\right)\right)=    
	\frac{1 }{m}\, \sum_{i=1}^{m-1} \, \mathcal F\left( \sin \left(\frac{\pi i  }{2m}\right)\right),
	\qquad m\in \NN,
\end{equation*}
so that we can use \cite[Proposition~13, p.~114]{A-D-M-Y} by which
\begin{equation*} 
	2 \widehat{\mathcal{S}}_{2m,1} = 1- 2\log(2) + \mathcal R\left(\frac{1}{2m}\right),
	\qquad m \in \NN,
\end{equation*}
and, thus, \eqref{limitatzeros1} holds for this case as well.

\textbf{2.2.2}
If  $k$ is odd then we proceed as in 1.2.2 and use \eqref{senoimpar}  in \eqref{difference11} to conclude that \eqref{limitatzeros1} holds  as well.

The proof of Theorem~\ref{teorema compatibilidad Ch-1} is complete.
\endproof

\emph{Proof of Proposition~\ref{prop1}.}

Let $ \theta/\pi = s/k $, 
with $s, k \in \NN, \; s<k$, and $ \operatorname{GCD} (s,k)=1 $. By Theorem~\ref{teorema compatibilidad Ch-1},
$$
{\mathcal D}_\infty(x) =  \log (2) + 2\,  \widehat{\mathcal{S}}_{k,s}.
$$
By (\ref{equality-ks-k1}), for Chebyshev polynomials of the first kind we can rewrite it as
$$
{\mathcal D}_\infty(x) =  \log (2) + 2\,  \widehat{\mathcal{S}}_{k,1}.
$$

If $k$ is even,  identity (\ref{identityR}) yields
$$
{\mathcal D}_\infty(x) =    1- \log (2)+\mathcal R\left(\frac{1}{k}\right) > 1- \log (2); $$
for the last inequality we have used  the fact that the coefficients in the series expansion (\ref{R-serie}) are all positive, so that 
 $\mathcal R(x)>0$  on $(0,1)$.

Analogously, when $k$ is odd, by  (\ref{newS}),
$${\mathcal D}_\infty(x) =  1- \log (2)+ 2\left [\mathcal R\left(\frac{1}{2k}\right)-\frac{1}{2}\mathcal R\left(\frac{1}{k}\right) \right] < 1- \log (2);
$$
now the inequality is the consequence of (\ref{convex function}), and the proposition is proved.

 \endproof

\end{document}